\documentclass{article}
\usepackage{amssymb,latexsym,amsmath,amsthm,mathrsfs}
\usepackage{graphicx}
\newcommand{\comment}[1]{}

\newcommand{\Int}{{\textstyle \int}}
\begin{document}
\title{A demonstration of a theorem on the order observed in the
sums of divisors\footnote{Originally published as
{\em Demonstratio theorematis circa ordinem in summis divisorum observatum},
Novi Commentarii Academiae scientiarum Imperialis Petropolitanae
\textbf{5} (1760), 75--83.
E244 in the Enestr{\"o}m index.
Translated from the Latin by Jordan Bell,
Department of Mathematics, University of Toronto, Toronto, Ontario, Canada.
Email: jordan.bell@gmail.com}}
\author{Leonhard Euler}
\date{}
\maketitle

\begin{center}
{\Large Summarium}
\end{center}
{\small Here the Celebrated Author shows completely what was still
desired in the preceding dissertation, and presently expounds a rigid
demonstration of this remarkable identity. Which, even if it rests on
common principles, still 
seems to display no small amount of ingenuity.
What follows from this argument has already been explained well enough
above.\footnote{Translator: The preceding paper in this volume of the
{\em Novi Commentarii} is E243, ``Observatio de summis divisorum'', where Euler states but does not prove
the pentagonal number theorem, and where he then derives the
recurrence for the sum of divisors function assuming it.}}

A while ago\footnote{Translator: I don't know how much later than E243 this paper
was written.} I discovered a theorem by which the nature of numbers has
been seen to be illuminated by no small amount, since in it is contained
an order which the sums of divisors arising from the numbers proceeding
in their natural series hold to each other.
For I showed that if all the divisors of each of the natural numbers
$1,2,3,4,5,6,7,8$ etc. are collected into one sum and these sums
of divisors are arranged in a series, which will be
\[
1,3,4,7,6,12,8,15,13,18,12,28,14,24,24,31,18 \quad \textrm{etc.},
\]
this series will be recurrent: each of the terms is determined from
the preceding according to a certain scale of relation. And this order
has not only been found to be highly noteworthy since scarcely anyone would
have suspected that this series would be bound by any fixed law, 
but also because at that time I was unable to discover any firm demonstration
of this order, even though I attempted it in many ways. I was led to find this
order when I was contemplating the expansion of the following infinite product
\[
s=(1-x)(1-x^2)(1-x^3)(1-x^4)(1-x^5)(1-x^6)(1-x^7) \; \textrm{etc.},
\]
and expanding this I concluded by induction that
\[
s=1-x-x^2+x^5+x^7-x^{12}-x^{15}+x^{22}+x^{26}-x^{35}-x^{40}+\textrm{etc.}
\]
The order of the exponents of $x$ is apparent by taking their differences;
the series of differences will be
\[
1, 1, 3, 2, 5, 3, 7, 4, 9, 5, 11, 6, 13, 7, 15, 8 \quad \textrm{etc.}
\]
Picking the terms alternately, it is clear that this series is admixed
from the series of odd
numbers and from the series of all integral numbers.
But indeed, that according to this law
$s=1-x-x^2+x^5+x^7-x^{12}-x^{15}+\textrm{etc.}$ if 
$s=(1-x)(1-x^2)(1-x^3)(1-x^4)(1-x^5)$ etc. I was only able to
show by induction, and I was not able to
show the equality with a solid demonstration.  
It was for this reason that I was not able to firmly demonstrate the
order which I found in the sums of divisors, but I indicated that
its demonstration depends on the demonstration of the equality between
the two infinite formulas that were exhibited above. But since
I have now obtained this demonstration, then also the order found
in the sums of divisors is no longer counted as I had
judged then among those truths which
are recognized yet can still not be demonstrated, but now merits
a place among the rigidly demonstrated truths. 
So that no doubt of this can remain, I will state and demonstrate
all the propositions on which the demonstration of this truth depends.

\begin{center}
{\Large Proposition 1}
\end{center}
{\em If
\[
s=(1+\alpha)(1+\beta)(1+\gamma)(1+\delta)(1+\epsilon)(1+\zeta)(1+\eta) \; \textrm{etc.}
\]
then this product arising from infinitely many factors may be
converted into the following series}
\[
\begin{split}
&s=(1+\alpha)+\beta(1+\alpha)+\gamma(1+\alpha)(1+\beta)+\delta(1+\alpha)(1+\beta)(1+\gamma)\\
&+\epsilon(1+\alpha)(1+\beta)(1+\gamma)(1+\delta)+\zeta(1+\alpha)(1+\beta)(1+\gamma)(1+\delta)(1+\epsilon)+\textrm{etc.}
\end{split}
\]

\begin{center}
{\Large Demonstration}
\end{center}

Since the first term of the series is $(1+\alpha)$ and the second is
$=\beta(1+\alpha)$, the sum of the first and second is $=(1+\alpha)(1+\beta)$;
now if the third term $\gamma(1+\alpha)(1+\beta)$ is added this will
yield $(1+\alpha)(1+\beta)(1+\gamma)$; and let the fourth term, which
is $\delta(1+\alpha)(1+\beta)(1+\gamma)$, be added, and the sum will be
\[
=(1+\alpha)(1+\beta)(1+\gamma)(1+\delta).
\]
And thus by proceeding to infinity, the whole sum of the entire series,
that is of all its terms, will be brought to this product
\[
(1+\alpha)(1+\beta)(1+\gamma)(1+\delta)(1+\epsilon)(1+\zeta) \; \textrm{etc.}
\]
Therefore it is clear that if
\[
s=(1+\alpha)(1+\beta)(1+\gamma)(1+\delta)(1+\epsilon)(1+\zeta) \; \textrm{etc.}
\]
then on the other hand
\[
s=(1+\alpha)+\beta(1+\alpha)+\gamma(1+\alpha)(1+\beta)+\delta(1+\alpha)(1+\beta)(1+\gamma)+\textrm{etc.}
\]

\begin{center}
{\Large Proposition 2}
\end{center}
{\em If
\[
s=(1-x)(1-x^2)(1-x^3)(1-x^4)(1-x^5)(1-x^6)\; \textrm{etc.}
\]
then this product arising from infinitely many factors 
may be reduced to the series}
\[
s=1-x-x^2(1-x)-x^3(1-x)(1-x^2)-x^4(1-x)(1-x^2)(1-x^3)-\textrm{etc.}
\]

\begin{center}
{\Large Demonstration}
\end{center}

If the form $s=(1-x)(1-x^2)(1-x^3)(1-x^4)(1-x^5)$ etc. is compared with the
preceding form $s=(1+\alpha)(1+\beta)(1+\gamma)(1+\delta)(1+\epsilon)$ etc.,
it is apparent that
\[
\alpha=-x, \quad \beta=-x^2, \quad \gamma=-x^3, \quad \delta=-x^4,
\quad \epsilon=-x^5 \quad \textrm{etc.}
\]
Then the truth of the proposition will be apparent by
substituting these given values into the series which was found
equal to the product $s$, namely that
\[
s=1-x-x^2(1-x)-x^3(1-x)(1-x^2)-x^4(1-x)(1-x^2)(1-x^3)-\textrm{etc.}
\]

\begin{center}
{\Large Proposition 3}
\end{center}
{\em If}
\[
s=(1-x)(1-x^2)(1-x^3)(1-x^4)(1-x^5)(1-x^6)(1-x^7) \; \textrm{etc.}
\]
{\em by expanding this infinite product by multiplication and by
arranging the terms according to powers of $x$ it will be}
\[
s=1-x^1-x^2+x^5+x^7-x^{12}-x^{15}+x^{22}+x^{26}-x^{35}-x^{40}+x^{51}+x^{57}
-\textrm{etc.},
\]
{\em whose rule of formation is the very one which was explained above.}

\begin{center}
{\Large Demonstration}
\end{center}
Since 
\[
s=(1-x)(1-x^2)(1-x^3)(1-x^4)(1-x^5)(1-x^6)(1-x^7) \; \textrm{etc.},
\]
it will be 
\[
s=1-x-x^2(1-x)-x^3(1-x)(1-x^2)-x^4(1-x)(1-x^2)(1-x^3)-\textrm{etc.}
\]
Let us put
\[
s=1-x-Ax^2;
\]
then it will be
\[
A=1-x+x(1-x)(1-x^2)+x^2(1-x)(1-x^2)(1-x^3)+\textrm{etc.}
\]
Let us expand all the terms by just the factor $1-x$ and
let us arrange them in the following way
\[
A=\bigg\{
\begin{array}{lllll}
&-x&-x^2(1-x^2)&-x^3(1-x^2)(1-x^3)&-\textrm{etc.}\\
1&+x(1-x^2)&+x^2(1-x^2)(1-x^3)&+x^3(1-x^2)(1-x^3)(1-x^4)&+\textrm{etc.}
\end{array}
\]
and by collecting together the terms written below it will be
\[
A=1-x^3-x^5(1-x^2)-x^7(1-x^2)(1-x^3)-x^9(1-x^2)(1-x^3)(1-x^4)-\textrm{etc.}
\]
Let us put
\[
A=1-x^3-Bx^5;
\]
it will be
\[
B=1-x^2+x^2(1-x^2)(1-x^3)+x^4(1-x^2)(1-x^3)(1-x^4)+\textrm{etc.};
\]
in all of these terms we let $1-x^2$ be expanded, and it will become
\[
B=\bigg\{
\begin{array}{lllll}
&-x^2&-x^4(1-x^3)&-x^6(1-x^3)(1-x^4)&-\textrm{etc.}\\
1&+x^2(1-x^3)&+x^4(1-x^3)(1-x^4)&+x^6(1-x^3)(1-x^4)(1-x^5)&+\textrm{etc.}
\end{array}
\]
and by again collecting the terms written below, we will have
\[
B=1-x^5-x^8(1-x^3)-x^{11}(1-x^3)(1-x^4)-x^{14}(1-x^3)(1-x^4)(1-x^5)-\textrm{etc.}
\]
Let us put
\[
B=1-x^5-Cx^8;
\]
it will be
\[
C=1-x^3+x^3(1-x^3)(1-x^4)+x^6(1-x^3)(1-x^4)(1-x^5)+\textrm{etc.},
\]
where we may expand the factor $1-x^3$ in all the terms, so that it will become,
writing it as above,
\[
C=\bigg\{
\begin{array}{lllll}
&-x^3&-x^6(1-x^4)&-x^9(1-x^4)(1-x^5)&-\textrm{etc.},\\
1&+x^3(1-x^4)&+x^6(1-x^4)(1-x^5)&+x^9(1-x^4)(1-x^5)(1-x^6)&+\textrm{etc.},
\end{array}
\]
which we collect to get
\[
C=1-x^7-x^{11}(1-x^4)-x^{15}(1-x^4)(1-x^5)-x^{19}(1-x^4)(1-x^5)(1-x^6)-
\textrm{etc.}
\]
Let us put
\[
C=1-x^7-Dx^{11};
\]
it will be
\[
D=1-x^4+x^4(1-x^4)(1-x^5)+x^8(1-x^4)(1-x^5)(1-x^6)+\textrm{etc.},
\]
which turns into this form
\[
D=\bigg\{
\begin{array}{lllll}
&-x^4&-x^8(1-x^5)&-x^{12}(1-x^5)(1-x^6)&-\textrm{etc.}\\
1&+x^4(1-x^5)&+x^8(1-x^5)(1-x^6)&+x^{12}(1-x^5)(1-x^6)(1-x^7)&+\textrm{etc.},
\end{array}
\]
and thus it will be
\[
D=1-x^9-x^{14}(1-x^5)-x^{19}(1-x^5)(1-x^6)-x^{24}(1-x^5)(1-x^6)(1-x^7)-\textrm{etc.}
\]
And if one puts next
\[
D=1-x^9-Ex^{14},
\]
it will similarly turn out
\[
E=1-x^{11}-Fx^{17}
\]
and then further
\[
F=1-x^{13}-Gx^{20}, \quad G=1-x^{15}-Hx^{23}, \quad H=1-x^{17}-Ix^{26}
\quad \textrm{etc.}
\]
Let us now successively replace these values, and it will be
\[
\begin{array}{rcr}
s&=&1-x-Ax^2,\\
Ax^2&=&x^2(1-x^3)-Bx^7,\\
Bx^7&=&x^7(1-x^5)-Cx^{15},\\
Cx^{15}&=&x^{15}(1-x^7)-Dx^{26},\\
Dx^{26}&=&x^{26}(1-x^9)-Ex^{40}\\
&\textrm{etc.}&
\end{array}
\]
From which we will have
\[
s=1-x-x^2(1-x^3)+x^7(1-x^5)-x^{15}(1-x^7)+x^{26}(1-x^9)-x^{40}(1-x^{11})+\textrm{etc.}
\]
or the very thing which is to be demonstrated,
\[
s=1-x-x^2+x^5+x^7-x^{12}-x^{15}+x^{22}+x^{26}-x^{35}-x^{40}+x^{51}+\textrm{etc.},
\]
from which at once the law of the exponents indicated above by differences
is clearly evident.

\begin{center}
{\Large Proposition 4,\\or,\\the Principal Theorem to be Demonstrated}
\end{center}
{\em If the notation $\Int n$ denotes the sum of all the divisors
of the number $n$ and similarly for all lesser numbers, just as for $n-\alpha$
it will be designated by $\Int (n-\alpha)$, then the sum of all the divisors
of $n$, or $\Int n$, will depend on the sums of the divisors of lesser numbers,
as}
\[
\begin{split}
&\Int n=\Int(n-1)+\Int(n-2)-\Int(n-5)-\Int(n-7)+\Int(n-12)+\Int(n-15)\\
&-\Int(n-22)-\Int(n-26)+\Int(n-35)+\Int(n-40)-\Int(n-51)-\Int(n-57)+\textrm{etc.}
\end{split}
\]

The following should be noted here:

1. The signs $+$ and $-$ alternately affect pairs of terms of this progression.

2. The law of the numbers $1,2,5,7,12,15,22,26$ etc. is clear from their
differences, which are $1,3,2,5,3,7,4$ etc.; from this one gathers that all
the terms are contained in the general formula $\frac{3zz \pm z}{2}$.

3. In each case, those terms of the progression are taken which remain
positive after the $\Int$ sign; while the others, for which the $\Int$ sign
comes in front of negative numbers, are to be omitted; thus for
$n=10$, it will be $\Int 10=\Int 9+\Int 8-\Int 5-\Int 3=13+15-6-4=18$.

4. The term $\Int(n-n)$ will occur in those cases in which $n$ is a number
from the series $1,2,5,7,12,15$ etc., and in these cases the given number
$n$ itself ought to be taken for the value of the term $\Int(n-n)$ or
$\Int 0$;
thus if $n=7$, it will be $\Int 7=\Int 6+\Int 5-\Int 2-\Int 0=12+6-3-7=8$,
and if $n=12$ then it will be $\Int 12=\Int 11+\Int 10-\Int 7-\Int 5+\Int 0=
12+18-8-6+12=28$.

\begin{center}
{\Large Demonstration}
\end{center}

Let us form the series
\[
z=x\Int 1+x^2\Int 2+x^3\Int 3+x^4\Int 4+x^5\Int 5+\textrm{etc.},
\]
where each power of $x$ is multiplied by the sum of the divisors
of the exponent of that power. 
Now if all the sums of divisors are resolved, it is clear that
the series will be transformed into this form
\[
\begin{array}{rll}
z=&1(x+x^2+x^3+x^4+x^5+\textrm{etc.})&+2(x^2+x^4+x^6+x^8+x^{10}+\textrm{etc.})\\
&+3(x^3+x^6+x^9+x^{12}+x^{15}+\textrm{etc.})&+4(x^4+x^8+x^{12}+x^{16}+x^{20}+\textrm{etc.})\\
&+5(x^5+x^{10}+x^{15}+x^{20}+x^{25}+\textrm{etc.})&+6(x^6+x^{12}+x^{18}+x^{24}+x^{30}+\textrm{etc.})\\
&&\textrm{etc.}
\end{array}
\]
By summing these geometric series it will become
\[
z=\frac{1x}{1-x}+\frac{2xx}{1-xx}+\frac{3x^3}{1-x^3}+\frac{4x^4}{1-x^4}
+\frac{5x^5}{1-x^5}+\frac{6x^6}{1-x^6}+\textrm{etc.}
\]
Let us multiply this form by $-\frac{dx}{x}$, and the
integral of the product will be
\[
-\int \frac{zdx}{x}=l(1-x)+l(1-xx)+l(1-x^3)+l(1-x^4)+l(1-x^5)+\textrm{etc.}
\]
or
\[
-\int\frac{zdx}{x}=l(1-x)(1-xx)(1-x^3)(1-x^4)(1-x^5)(1-x^6) \; \textrm{etc.};
\]
but as the expression after the logarithm sign is the same as that
in the preceding proposition 
called $=s$, it will be $-\int \frac{zdx}{x}=ls$, and hence by taking
the other value for $s$ it will also be
\[
-\int \frac{zdx}{x}=l(1-x-x^2+x^5+x^7-x^{12}-x^{15}+x^{22}+x^{26}-\textrm{etc.}),
\]
whose differential divided by $\frac{-dx}{x}$ gives another value for $z$,
namely
\[
z=\frac{1x+2x^2-5x^5-7x^7+12x^{12}+15x^{15}-22x^{22}-\textrm{etc.}}{1-x-x^2+x^5+x^7-x^{12}-x^{15}+x^{22}+\textrm{etc.}};
\]
if this value is put equal to assumed one and both sides are multiplied by
the denominator $1-x-x^2+x^5+x^7-x^{12}$ etc., 
with the terms arranged according to powers of $x$ and collecting all on
one side it will turn out that:

{\scriptsize
\[
\begin{array}{rrrrrrrrrrrr}
0=&x\Int 1&+x^2\Int 2&+x^3\Int 3&+x^4\Int 4&+x^5\Int 5&+x^6\Int 6&+x^7\Int 7
&+x^8\Int 8&+x^9\Int 9&+x^{10}\Int 10&+\textrm{etc.}\\
&&-\Int 1&-\Int 2&-\Int 3&-\Int 4&-\Int 5&-\Int 6&-\Int 7&-\Int 8&-\Int 9&\\
&&&-\Int 1&-\Int 2&-\Int 3&-\Int 4&-\Int 5&-\Int 6&-\Int 7&-\Int 8&\\
&&&&&&+\Int 1&+\Int 2&+\Int 3&+\Int 4&+\Int 5&\\
&&&&&&&&+\Int 1&+\Int 2&+\Int 3&\\
&&&&&&&&&&\vdots&\\
&-1&-2&*&*&+5&*&+7&*&*&*&
\end{array}
\]
}
Whence, with the coefficients of all the powers of $x$ equal to $0$,
it follows that
\[
\begin{array}{ll}
\Int 1=1,&\Int 6=\Int 5+\Int 4-\Int 1,\\
\Int 2=\Int 1+2,&\Int 7=\Int 6+\Int 5-\Int 2-7,\\
\Int 3=\Int 2+\Int 1,&\Int 8=\Int 7+\Int 6-\Int 3-\Int 1,\\
\Int 4=\Int 3+\Int 2,&\Int 9=\Int 8+\Int 7-\Int 4-\Int 2,\\
\Int 5=\Int 4+\Int 3-5,&\Int 10=\Int 9+\Int 8-\Int 5-\Int 3;
\end{array}
\]
and by a light inspection the character of this equation will clearly be
\[
\Int n=\Int(n-1)+\Int(n-2)-\Int(n-5)-\Int(n-7)+\Int(n-12)+\Int(n-15)-\textrm{etc.}
\]
In each case this progression is continued until it comes to sums of negative
numbers. Then it is clear by itself that the actual numbers $1,2,5,7$ etc.
which are seen in these formulas take the place of the term $\Int 0$;
from which one concludes that in the cases in which the term
$\Int(n-n)$ or $\Int 0$ occurs in the progression found for
$\Int n$, its value is always to be taken equal to the given number
$n$ itself; and thus a complete and perfect demonstration of
the proposed theorem is obtained, which, since beyond the treatment
of infinite series it proceeds by logarithms and differentials, is indeed
less natural, but because of this it is to be considered all the more notable.

\end{document}